\documentclass[
a4paper,
11pt,
DIV=14,
toc=bibliography, 
headsepline, 
parskip=half-,
abstract=true,
]{scrartcl}

\usepackage[centertags]{amsmath}
\usepackage[utf8]{inputenc}
\usepackage{enumerate,amsfonts,amssymb,amsthm,amsopn,cite,array,comment}

\usepackage[usenames]{color}
\definecolor{citegreen}{rgb}{0,0.6,0}
\definecolor{refred}{rgb}{0.8,0,0}
\usepackage[colorlinks, citecolor=citegreen, linkcolor=refred, urlcolor=black]{hyperref}

\title{Perelman's entropy and Gromov--Hausdorff versus dynamical stability of round spheres}
\author{Reto Buzano and Shengwen Wang}
\date{}

\newtheorem{thm}{Theorem}[section]
\newtheorem{lem}[thm]{Lemma}

\newtheorem{prop}[thm]{Proposition}

\theoremstyle{definition}
\newtheorem{rem}[thm]{Remark}

\newtheorem{defn}[thm]{Definition}

\numberwithin{equation}{section}

\newcommand{\eps}{\varepsilon}
\newcommand{\R}{\mathbb{R}}
\newcommand{\N}{\mathbb{N}}
\newcommand{\Sph}{\mathbb{S}}
\newcommand{\Nash}{\mathcal{N}}
\newcommand{\W}{\mathcal{W}}
\newcommand{\Vol}{\operatorname{Vol}}
\newcommand{\diam}{\operatorname{diam}}
\newcommand{\Area}{\operatorname{Area}}
\newcommand{\Ric}{\operatorname{Ric}}
\newcommand{\Rm}{\operatorname{Rm}}
\newcommand{\dGH}{d_{\mathrm{GH}}}
\providecommand{\abs}[1]{\lvert #1\rvert}
\providecommand{\norm}[1]{\lVert #1\rVert}

\newcommand\printaddress{{
\setlength{\parindent}{15pt}
\footnotesize~
\par
{\scshape Reto Buzano}
\newline 
Universit\`a di Torino, 
Dipartimento di Matematica,
Via Carlo Alberto 10, 
10123 Torino, Italy 
\newline
\emph{E-mail address:} 
\texttt{reto.buzano@unito.it}
\par
{\scshape Shengwen Wang}
\newline 
Queen Mary University of London, 
School of Mathematical Sciences, 
Mile End Road, 
London E1 4NS, UK
\newline
\textit{E-mail address:} 
\texttt{shengwen.wang@qmul.ac.uk}
\par
}}

\begin{document}
\maketitle
\begin{abstract}
We study the Gromov--Hausdorff stability of round spheres under perturbations of Perelman's $\nu$-entropy. We show that in dimension two almost maximal entropy, together with an area normalisation, implies Gromov--Hausdorff closeness to the round two-sphere. In contrast, in dimensions $n \geq 3$, we construct smooth metrics on $\Sph^n$ whose entropies converge to that of the round sphere, but which remain a definite Gromov--Hausdorff distance away from it. These metrics can be chosen conformal and uniformly bi-Lipschitz equivalent to the round metric, with the same volume as the round sphere. Finally, in dimension $n=3$, we contrast this Gromov--Hausdorff instability with a dynamical stability result: metrics with $\nu$-entropy above that of the round shrinking cylinder converge exponentially fast to round metrics under the volume-normalised Ricci flow.
\end{abstract}


\section{Introduction and main results}

Perelman's entropy functionals play a central role in the study of Ricci flow, in particular in the analysis of singularities and gradient shrinking Ricci solitons \cite{Per}. Besides their importance for the flow, these functionals also provide important geometric invariants of Riemannian metrics. This raises the natural question of how much information about the (metric) geometry of a manifold can be recovered from the entropy functionals. In particular, one may ask whether metrics whose entropy is close to that of a distinguished shrinking Ricci soliton must themselves be geometrically close to that soliton. 

Here, we work with Perelman's $\nu$-entropy functional. In \cite[Conjecture 9.3]{BWa}, Wang conjectured that in all dimensions $n$, this entropy functional attains its maximal value among all smooth metrics on closed manifolds at the round $n$-sphere. He further conjectured that if the entropy of a smooth closed Riemannian manifold $(M,g)$ is sufficiently close to this maximal value, then $M$ is diffeomorphic to the round sphere and the normalised Ricci flow starting at $(M,g)$ converges to a round sphere. In \cite[Theorem 1.4]{BWa}, he proved such a dynamical stability result under a \emph{Ricci curvature lower bound} and volume assumptions, obtaining in particular a quantitative distance comparison between the distance induced by $g$ and the one induced by the round metric, thus also obtaining \emph{stability in the Gromov--Hausdorff sense}. Immediately before Conjecture 9.3, Wang suggests that the geometric hypotheses in his theorems might eventually be replaced by entropy conditions alone. 

In this article, we first study the Gromov--Hausdorff stability of round spheres under such an entropy bound and a volume normalisation -- as $\nu$ is scale-invariant, any question concerning Gromov--Hausdorff stability must include a normalisation of scale. We show that, once a volume normalisation is imposed, the Gromov--Hausdorff stability question exhibits a striking dimensional dichotomy: round spheres are Gromov--Hausdorff stable under small perturbations of $\nu$ in dimension two (Theorem \ref{thm.main1}), but this stability fails in every higher dimension $n \geq 3$, even among uniformly bi-Lipschitz conformal perturbations (Theorem \ref{thm.main2}). This demonstrates that the quantitative distance comparison between initial metrics $g$ and round metrics, obtained in Wang's Theorem 1.4 under Ricci curvature bounds, cannot be recovered from near-spherical $\nu$-entropy and volume normalisation alone in dimensions $n\geq 3$.

To state our main results precisely, we denote by $g_0$ the unit round metric on $\Sph^n$, and let $\nu_n=\nu(g_0)$. We first state our stability result in dimension two.

\begin{thm}[Gromov--Hausdorff stability in dimension $n=2$]\label{thm.main1}
Let $(M_k,g_k)_{k \in \N}$ be a sequence of smooth, closed Riemannian surfaces satisfying $\nu(g_k)\to \nu_2$ and
\begin{equation*}
	\Area(M_k,g_k)=4\pi.
\end{equation*}
Then, for all sufficiently large $k$, the surfaces $M_k$ are diffeomorphic to $\Sph^2$, and
\begin{equation*}
	\dGH \big((M_k,g_k), (\Sph^2,g_0)\big)\to 0.
\end{equation*}
\end{thm}

This stability result is contrasted by the following instability result in higher dimensions.

\begin{thm}[Gromov--Hausdorff instability in dimension $n \geq 3$]\label{thm.main2}
For every fixed $n \geq 3$, there exists a sequence $(g_k)_{k\in \N}$ of smooth Riemannian metrics on $\Sph^n$ with $\nu(g_k)\to \nu_n$, but
\begin{equation*}
	\liminf_{k\to\infty} \dGH \big((\Sph^n,g_k),(\Sph^n,g_0)\big)>0.
\end{equation*}
Moreover, the metrics $g_k$ can be chosen conformal to $g_0$, uniformly bi-Lipschitz equivalent to $g_0$, and with volumes equal to that of the round metric $g_0$.
\end{thm}

It is interesting to compare these results with the corresponding \emph{extrinsic} stability question for the entropy introduced by Colding and Minicozzi in mean curvature flow \cite{CM}. Bernstein and Wang \cite{BW} established Hausdorff stability of the round two-sphere under small perturbations of that entropy. The second author then extended their stability result to round hyperspheres in all higher dimensions \cite{SWa}. This comparison highlights an interesting contrast between the extrinsic Colding--Minicozzi entropy, for which Hausdorff stability \emph{holds in all dimensions}, and Perelman's intrinsic $\nu$-entropy, for which Gromov--Hausdorff stability \emph{fails in dimensions at least three}.

We point out that the instability in Theorem \ref{thm.main2} persists above the entropy of the round shrinking cylinder. Indeed, the entropies of the round spheres satisfy
\begin{equation}
	\nu_n>\nu_{n-1} =\nu_c(\Sph^{n-1}\times\R), \qquad n \geq 3.
\end{equation}
Our metrics therefore eventually have entropy strictly greater than that of the corresponding round cylinder, but still remain a definite Gromov--Hausdorff distance away from the round sphere.

Despite our Gromov--Hausdorff instability result and the failure of controlling the initial geometry via only an entropy bound, we are still able to prove a dynamical stability result and give an affirmative answer to Wang's Conjecture 9.3 in dimension $n=3$, based on Perelman's singularity theory and the classification of singularity models. More precisely, we prove the following.

\begin{thm}[Ricci flow stability in dimension $n = 3$]\label{thm.main3}
Among all smooth, closed Riemannian three-manifolds $(M, g)$ we have $\nu(M,g) \leq \nu_3$ with equality if and only if $(M,g)$ is isometric to the round sphere $(\Sph^3,g_0)$, up to scaling. Furthermore, if $\nu(M,g) > \nu_2 = \nu_c(\Sph^2 \times \R)$, then $M$ is diffeomorphic to $\Sph^3$ and the volume-normalised Ricci flow starting at $(M,g)$ exists for all positive times and converges exponentially fast to a round metric.
\end{thm}

An analogous result to Theorem \ref{thm.main3} also holds in dimension $n=2$ by the topological part of Theorem \ref{thm.main1} and the Hamilton-Chow long-time existence and convergence theorem for normalised Ricci flow on topological two-spheres \cite{Ham88, Cho}.

Theorems \ref{thm.main2} and \ref{thm.main3} reveal a sharp distinction between geometric and dynamical stability under near-maximal entropy in dimension three. While such an entropy bound does not force the initial metric to be close to the round sphere in the Gromov--Hausdorff topology, it nevertheless guarantees convergence to a round metric under the normalised Ricci flow. In particular, the entropy gap to the round shrinking cylinder suffices to exclude non-round singularity formation, despite the absence of any quantitative control on the initial distance geometry. Moreover, the entropy gap appearing in Theorem \ref{thm.main3} is explicit:
\begin{equation}
	\nu_3-\nu_2=\frac{\log\pi-1}{2} \approx 0.07236 >0.
\end{equation}

\subsection{Comparison to previous Gromov--Hausdorff (in)stability results}

Classical stability results in Riemannian geometry often assert that almost equality of a geometric quantity or in a geometric comparison theorem forces closeness to a model space in the Gromov--Hausdorff topology. A fundamental example is Colding's almost-rigidity theorem \cite{Col}, which shows that a closed manifold satisfying $\Ric \geq n-1$ and having volume sufficiently close to that of the unit round sphere must be Gromov--Hausdorff close to the round sphere. 

The situation becomes considerably more delicate when one replaces lower Ricci curvature bounds by weaker geometric or analytic assumptions, such as lower scalar curvature bounds or, as in our case, bounds on Perelman's entropy. In particular, it turns out that geometric features concentrated in regions of arbitrarily small volume can have a substantial effect on distances, even when their effect on certain integral or analytic quantities is negligible.

The results most closely related to ours are those of Lee, Naber and Neumayer \cite{LNN}, who investigated geometric stability under almost Euclidean lower bounds for the scalar curvature and Perelman's entropy. They showed that such bounds do not, in general, imply Gromov--Hausdorff closeness to Euclidean space. In particular, they constructed examples in dimensions $n\geq4$ for which the scalar curvature and entropy restricted to scales $0<\tau<2$ are almost nonnegative, while the associated distance functions exhibit substantial degeneration. On the other hand, they also introduced a weaker notion of convergence, called $d_p$-convergence, and established corresponding $\eps$-regularity and compactness theorems. Our examples from Theorem \ref{thm.main2} exhibit a related failure of Gromov--Hausdorff stability, but in a different setting. We consider metrics on the compact manifold $\Sph^n$, and the relevant quantity is the global $\nu$-entropy approaching the value of the round shrinking sphere, rather than a scale-restricted entropy approaching its Euclidean value. Moreover, our examples are conformal to the round metric and remain uniformly bi-Lipschitz equivalent to it. Thus, the failure of Gromov--Hausdorff convergence occurs without collapse or degeneration of the underlying metric topology, and already in dimension $n=3$.

Another closely related result appears in the work of Kazaras and Xu \cite{KX}, who constructed so-called \emph{drawstring metrics} along codimension-two submanifolds while preserving scalar curvature lower bounds up to an arbitrarily small error. In particular, their construction allows distances along a prescribed submanifold to be shortened by a fixed factor, producing nonstandard Gromov--Hausdorff limits. For a fixed finite shortening factor, their approximating metrics can also be chosen uniformly bi-Lipschitz equivalent to the background metric. The main differences between their results and ours lie in the precise construction and the analytic quantity being controlled: their metrics involve different deformations in the tangential and normal directions of the submanifold while our metrics are obtained by a conformal deformation; their construction almost preserves a scalar curvature bound while our construction almost preserves Perelman's $\nu$-entropy.

Several recent results illustrate the importance of the precise hypotheses and topology in geometric (in)stability results. Neumayer \cite{Neu} showed that Gromov--Hausdorff closeness to Euclidean space can be recovered from almost Euclidean scalar curvature and entropy bounds if one additionally assumes sufficiently sharp upper bounds on the volumes of geodesic balls. Hirsch and Zhang \cite{HZ} established quantitative stability of Llarull's scalar curvature rigidity theorem, obtaining closeness to spherical geometry outside suitably controlled exceptional sets. Examples of Sweeney \cite{Swe} demonstrate that Gromov--Hausdorff convergence may fail in scalar curvature stability problems where weaker notions of convergence remain valid. Very interesting results in a slightly different direction were also obtained by Lee and Topping \cite{LT}, who showed that, in dimensions at least four, arbitrary metrics conformal to the round sphere can arise as uniform distance limits of metrics with positive scalar curvature. 

\subsection{Outline of the proofs and organisation of the article}

Our construction of metrics $g_k$ for Theorem \ref{thm.main2} is based on modifying the round metric on smaller and smaller tubular neighbourhoods of a totally geodesic codimension two sphere $\Sigma^{n-2} \subset\Sph^n$. More precisely, we construct conformal metrics $g_k=e^{2 \varphi_k}g_0$ that agree with the round metric outside small tubular neighbourhoods of $\Sigma$, while restricting to a fixed multiple of the round metric along $\Sigma$. Consequently, distances between points of $\Sigma$ remain uniformly shortened, even as the region in which the metric is modified shrinks to $\Sigma$. This produces a nontrivial distortion of the limiting distance function.

The main analytic difficulty is to show that these deformations have an asymptotically negligible effect on Perelman's $\nu$-entropy. The essential observation is that the conformal factors can be chosen so that their gradients are small in a suitable sense, despite the fact that the conformal factors have a fixed nonzero oscillation. This is achieved using \emph{double-logarithmic cutoff functions} and a \emph{critical Hardy inequality} in the two-dimensional normal directions to $\Sigma$, which is where the codimension-two condition comes in. In particular, we obtain estimates of the form
\begin{equation}
	\int_{\Sph^n} \abs{\nabla \varphi_k}^2 \psi^2\,dV_{g_0} \leq o(1)\int_{\Sph^n} \big(\abs{\nabla \psi}^2+\psi^2\big)\,dV_{g_0},
\end{equation}
for smooth functions $\psi$, uniformly in $k$, see Proposition \ref{prop.gradient}. These estimates allow us to compare the $\W$-functionals of $g_k$ and $g_0$ and ultimately show that their $\nu$-entropies converge. Since the $\nu$-functional involves an infimum over all positive scales, an important part of the argument is to obtain uniform control over the scales relevant to this infimum.

The preceding mechanism does not produce a corresponding instability in dimension two: in this case, a codimension-two submanifold is zero-dimensional, hence we no longer obtain a positive-dimensional submanifold along which distances are shortened. More importantly, in dimension two, Perelman's entropy provides a lower bound for the first eigenvalue of the Schrödinger type operator $-4\Delta_g+R_g$, which allows us to exploit spectral almost-rigidity results. More precisely, testing the $\W$-functional with a normalised first eigenfunction, and using the area normalisation, shows that entropy close to $\nu_2$ forces this first eigenvalue to be close to $2$. After an appropriate rescaling, the spectral almost-rigidity theorem of Xu \cite{Xu} yields the desired Gromov–Hausdorff convergence. The topological statement follows from a simple Gauss-Bonnet argument.

Finally, for the dynamical stability result in dimension three, we use Perelman's singularity theory \cite{Per}. The key observation is that the entropy bound strictly above $\nu_2$ on the one hand forces a finite time singularity for the un-normalised Ricci flow and on the other hand excludes noncompact singularity
models. As compact singularity models have strictly positive sectional curvatures, we can deduce that our Ricci flow will eventually have strictly positive Ricci curvature. At this point, the result follows by Hamilton's convergence theorem for three-manifolds with positive Ricci curvature \cite{Ham82}.

The paper is organised as follows. In Section \ref{sec.preliminaries}, we recall the relevant properties of Perelman's entropy and its values for round spheres, cylinders, and quotients. In Section \ref{sec.instability}, we first construct the conformal metrics used in the higher-dimensional counterexamples and establish the logarithmic Hardy estimate that controls their entropy. We then prove that the resulting metrics remain Gromov--Hausdorff separated from the round sphere, proving Theorem \ref{thm.main2}. In Section \ref{sec.stability}, we establish the stability theorem in dimension $n=2$, Theorem \ref{thm.main1}, using the spectral almost-rigidity result of Xu \cite{Xu}. Finally, in Section \ref{sec.dynamical}, we give a proof of the dynamical stability result from Theorem \ref{thm.main3} in dimension $n=3$.

\paragraph{AI disclosure statement.} The main ideas, mathematical approach, and final writing of this paper are due to the authors. During preparation, the authors used publicly available LLMs to assist with proofreading, literature searches, and the simplification of some estimates, particularly in the proof of Theorem \ref{thm.main2}. All AI-generated suggestions incorporated into the manuscript were verified and revised by the authors, who take full responsibility for the content of this publication.

\section{Preliminaries and notation}\label{sec.preliminaries}

Let $M$ be a smooth, closed manifold. (We always assume that manifolds are connected.) We denote by $g$ a smooth Riemannian metric on $M$, by $f$ a smooth real-valued function on $M$, and by $\tau>0$ a scaling parameter. In his ground-breaking article \cite{Per}, Perelman introduced the $\W$-functional
\begin{equation*}
	\W(g,f,\tau) =\int_M \left[\tau(R_g+\abs{\nabla f}_g^2)+f-n\right] (4\pi\tau)^{-n/2}e^{-f} \,dV_g,
\end{equation*}
restricted to functions $f$ satisfying the normalisation
\begin{equation*}
	\int_M (4\pi\tau)^{-n/2}e^{-f}\,dV_g=1.
\end{equation*}
In the following, we often write $u = (4\pi\tau)^{-n/2}e^{-f} > 0$ and $v=\sqrt{u}$. This yields
\begin{equation*}
	f= -2\log  v-\frac n2\log (4\pi\tau),\qquad  u \abs{\nabla f}_g^2=4 \abs{\nabla v}_g^2,
\end{equation*}
and hence, in terms of $v > 0$, we can write
\begin{equation}\label{eq.W-v}
	\W(g,v,\tau)=\tau\int_M\left( 4\abs{\nabla v}_g^2+R_gv^2\right) dV_g
	-\int_M v^2\log  v^2\,dV_g-n-\frac{n}{2} \log (4\pi\tau),
\end{equation}
with normalisation
\begin{equation*}
	\norm{v}^2_{L^2(g)} = \int_M v^2\,dV_g=1.
\end{equation*}
For convenience, we will often decompose $\W$ as 
\begin{equation}\label{eq.W-decomp1}
	\W_g(v,\tau) = \W(g,v,\tau) =\tau Q_g(v)-\Nash_g(v^2)-n-\frac{n}{2} \log (4\pi\tau)
\end{equation}
where 
\begin{equation}\label{eq.W-decomp2}
	Q_g(v)=\int_M (4\abs{\nabla v}_g^2+R_gv^2)dV_g,\qquad  \Nash_g(v^2)=\int_M v^2\log  v^2dV_g.
\end{equation}
\begin{defn}
The \emph{$\mu$-entropy at scale $\tau$} and the scale-invariant \emph{$\nu$-entropy} are defined by
\begin{equation}\label{eq.entropies}
	\mu(g,\tau) =\inf_v \W(g,v,\tau), \qquad
	\nu(g) =\inf_{\tau>0}\mu(g,\tau),
\end{equation}
where the infimum in the definition of $\mu(g,\tau)$ is taken over all strictly positive $v\in C^\infty(M)$ with $\norm{v}_{L^2(g)}^2=1$.
\end{defn}
\begin{rem}
Note that in the definition of $\mu(g,\tau)$ we can also allow properly normalised $v \in W^{1,2}(M)$ with $v \geq 0$ (using the convention that $0 \log 0 = 0$). This is due to density of smooth and strictly positive functions in this Sobolev space, as well as continuity of the $\W$-functional under strong $W^{1,2}$-convergence. By standard elliptic regularity and positivity results for the Euler-Lagrange equation, every nonnegative minimiser is smooth and strictly positive.
\end{rem}
Perelman's $\nu$-entropy is a fundamental tool in his proof of the Poincaré conjecture and, more generally, in modern Ricci flow theory, due to the fact that it is nondecreasing along Ricci flow. Moreover, when it is finite, equality on a nontrivial time interval characterises gradient shrinking Ricci solitons. In this article, we are mainly interested in this quantity as a geometric invariant for a fixed, non-evolving Riemannian metric $g$. Here, we briefly collect some basic properties of the $\nu$-functional in this setting.

Let $(\Sph^n,g_0)$ denote the unit round sphere, with induced distance function $d_{g_0}$. Its scalar curvature is $R_{g_0}=n(n-1)$, and we write $V_n=\Vol(\Sph^n,g_0)$. Its $\nu$-entropy (in fact, more generally, the $\nu$-functional of any positive Einstein manifold) is realised by a properly normalised constant function $\bar{v}$. For $\bar{v}=V_n^{-1/2}$, we find by \eqref{eq.W-decomp1} and \eqref{eq.W-decomp2},
\begin{equation*}
\mathcal W_{g_0}(\bar{v},\tau) =n(n-1)\tau+\log V_n-n-\frac n2\log(4\pi\tau).
\end{equation*}
This quantity is minimised for $\bar{\tau}=\frac1{2(n-1)}$ with value
\begin{equation}\label{eq.sphere-entropy}
    \nu_n:=\nu(g_0) =\log V_n +\frac n2\log\frac{n-1}{2\pi} -\frac{n}{2}.
\end{equation}
Moreover, the spherical entropies are strictly increasing with dimension and satisfy
\begin{equation}\label{eq.sphere-monotonicity}
    \nu_2=\log2-1 <\nu_3=\log(2\sqrt{\pi})-\frac{3}{2} <\cdots< \lim_{n\to\infty}\nu_n =\frac{\log2-1}{2}<0.
\end{equation}
These formulas and their properties are discussed in \cite{CHI}.

\begin{defn}
For a \emph{complete} Riemannian $n$-manifold $(M,g)$, we set
\begin{equation*}
    \nu_c(g) := \inf_{\tau>0} \inf_v \W_g(v,\tau),
\end{equation*}
where the second infimum is taken over all $v\in C_c^\infty(M)$ with $v\geq 0$ and $\int_M v^2\,dV_g=1$, and where the value of $v^2\log v^2$ appearing in \eqref{eq.W-v} is interpreted as $0$ wherever $v$ vanishes.
\end{defn}
We have $\nu_c=\nu$ on closed manifolds and $\nu_c$ is invariant under diffeomorphisms and rescalings of the metric with a constant. For the round shrinking cylinder $(\Sph^{2}\times\R, g_{\mathrm{cyl}})$, where $g_{\mathrm{cyl}}=g_0+dz^2$ is the canonical product metric (here $g_0$ denotes the unit round metric on $\Sph^2$), we have
\begin{equation}\label{eq.cylinder-entropy}
	\nu_c(\Sph^{2}\times\R, g_{\mathrm{cyl}}) \leq \nu_{2}.
\end{equation}
Let us quickly sketch the proof. Since $R_{g_{\mathrm{cyl}}}=2$ and $\Vol(\Sph^2,g_0)=4\pi$, at scale $\bar{\tau}=\frac{1}{2}$, consider the normalised Gaussian
\begin{equation*}
	v_{\mathrm{cyl}}(\theta,z)^2 =\frac{1}{4\pi\sqrt{2\pi}}e^{-z^2/2}.
\end{equation*}
Using \eqref{eq.W-v}--\eqref{eq.W-decomp2} also for noncompactly supported functions whenever the corresponding integrals converge, a direct computation gives
\begin{equation*}
	Q_{g_{\mathrm{cyl}}}(v_{\mathrm{cyl}})=3,\qquad
	\Nash_{g_{\mathrm{cyl}}}(v_{\mathrm{cyl}}^2) =-\log(4\pi)-\frac{1}{2}\log(2\pi)-\frac{1}{2},
\end{equation*}
and hence
\begin{equation*}
	\W_{g_{\mathrm{cyl}}}\big( v_{\mathrm{cyl}},\tfrac{1}{2}\big) =\log2-1=\nu_2.
\end{equation*}
Now choose a cutoff function $\eta\in C_c^\infty(\R)$ with $0\leq\eta\leq1$ and $\eta\equiv1$ near $0$, and set
\begin{equation*}
    v_L =\frac{\eta(z/L)v_{\mathrm{cyl}}}{\norm{\eta(z/L)v_{\mathrm{cyl}}}_{L^2(g_{\mathrm{cyl}})}}.
\end{equation*}
Gaussian decay gives
\begin{equation*}
	\W_{g_{\mathrm{cyl}}}\big( v_L,\tfrac{1}{2}\big) 
	\to \W_{g_{\mathrm{cyl}}}\big( v_{\mathrm{cyl}},\tfrac{1}{2}\big)
\end{equation*}
and since each $v_L$ is an admissible compactly supported test function, \eqref{eq.cylinder-entropy} follows.

In fact, the round shrinking cylinders equipped with their canonical product metrics have entropy
\begin{equation*}
    \nu_c(\Sph^{n-1}\times\R,g_{\mathrm{cyl}}) =\nu_{n-1}<\nu_n 
\end{equation*}
for all $n\geq 3$. An inequality follows in the same way as in the proof above for $n=3$. The reverse inequality follows from a sharp logarithmic Sobolev inequality for gradient shrinking Ricci solitons, see e.g.~\cite[Theorem 1.1]{LW20}.

Finally, if $\Gamma$ is a finite group acting freely and isometrically on the round sphere, then
\begin{equation}\label{eq.quotient-entropy}
    \nu(\Sph^n/\Gamma) =\nu_n-\log\abs{\Gamma}.
\end{equation}
This identity follows from the corresponding formulas for the normalised potentials of gradient shrinking Ricci solitons.

It is conjectured that the round $n$-spheres maximise the $\nu$-entropy among all smooth and closed $n$-dimensional manifolds for all dimensions $n$. For $n=2$ and $n=3$, we will see proofs of these results below. It is easy to see that in higher dimensions $n\geq 4$, the round sphere maximises $\nu$ among closed positive Einstein manifolds. Indeed, after rescaling an Einstein metric $g$ so that $\Ric_g=(n-1)g$, testing $\W$ with the constant normalised function $v=\Vol(M,g)^{-1/2}$ at $\bar{\tau}=\frac1{2(n-1)}$ gives
\begin{equation*}
    \nu(g)\leq \nu_n+ \log\frac{\Vol(M,g)}{V_n} \leq\nu_n,
\end{equation*}
where the last inequality follows from Bishop's volume comparison theorem. Equality is possible only for the round sphere. Moreover, the round metric is a \emph{local} maximum of the $\nu$-entropy among smooth metrics in every dimension. This follows from the stability of the round sphere under the normalised Ricci flow (see \cite{CHI}), together with Perelman's entropy monotonicity.

\section{Gromov--Hausdorff instability of $\Sph^n$ for $n\geq 3$}\label{sec.instability}

In this section we prove Theorem \ref{thm.main2}. As explained in the introduction, Lee--Naber--Neumayer \cite{LNN} and Kazaras--Xu \cite{KX} modify a Riemannian metric on a homogeneous manifold by shortening the distance near a submanifold of codimension three or two, respectively, while keeping it fixed in the rest of the manifold. Their main goal is to control the scalar curvature of the modified metric, but \cite{LNN} also proves an estimate of the $\mu$-entropy for a fixed scale $\tau$. Inspired by their work, we will shorten the round metric on $\Sph^n$ near a totally geodesic codimension two sphere $\Sigma^{n-2}$. Unlike \cite{LNN} and \cite{KX}, we use a conformal deformation, with the aim of controlling the scale-invariant $\nu$-entropy of the modified metric.

Fix $n\geq 3$ and let $\Sph^n$ be the standard unit sphere in $\R^{n+1}$ with its (induced) round metric $g_0$. We denote by $\Sigma$ a totally geodesic codimension two sphere in $\Sph^n$, for example we could simply pick 
\begin{equation*}
	\Sigma=\{ (x_1,\ldots,x_{n+1}) \in \Sph^n \mid x_n=x_{n+1}=0\}.
\end{equation*}

The following construction is motivated by the classical logarithmic cutoff trick in potential theory. Choose a smooth Heaviside type function $\chi: \R \to [0,1]$ with
\begin{equation*}
	\chi (s)=\begin{cases}0, \quad &s\in (-\infty, 0],\\
	1,\quad &s\in [1, \infty).\end{cases}
\end{equation*}
We also require $\chi'(s)\geq 0$ for all $s \in \R$ and fix a constant $0<c_n<-\frac{\nu_n}{n}$, for example $c_n=-\frac{\nu_n}{2n}$.

\begin{defn} 
We define the sequence of conformal metrics $g_k=e^{2\varphi_k}g_0$ on $\Sph^n$ by 
\begin{equation}\label{eq.conffactors}
	\varphi_k(p):=\begin{cases}
	-c_n,\quad &r=0,\\
	-c_n\cdot\chi\left(\displaystyle\frac{\log \log (1/r)-k}{k}\right),\quad &0<r<r_k=e^{-e^k},\\
	0, \quad &r_k=e^{-e^k}\leq r\leq\frac{\pi}{2}.
\end{cases}.
\end{equation}
\end{defn}
Here and in the following $r = r(p):=d_{g_0}(p,\Sigma)$ denotes the Riemannian distance in $\Sph^{n}$ to the submanifold $\Sigma$. Note that $\varphi_k$ are smooth (but not analytic) functions on $\Sph^n$ satisfying $-c_n\leq \varphi_k \leq0$, and the conformal metrics only differ from $g_0$ in a small $r_k$ tubular neighbourhood of the geodesic sphere $\Sigma$. On $\Sigma$, the metric is multiplied with the fixed factor $e^{-2c_n}<1$.

\begin{rem}
The condition $0<c_n<-\frac{\nu_n}{n}$ is not needed for the geometric construction above, the Gromov--Hausdorff distance from the round metric $g_0$ or the Hardy estimate below; it enters only when we want to exclude degenerating entropy scales in the proof of Proposition \ref{prop.conv-entropies}.
\end{rem}

We first prove the following key gradient estimate for $\varphi_k$.
\begin{prop}[Gradient estimates for the conformal factors]\label{prop.gradient}
There exists a constant $C=C(n,\chi)$ such that, for all sufficiently large $k$ and all functions $\psi \in C^\infty(\Sph^n)$,
\begin{equation*}
	\int_{\Sph^n} \abs{\nabla \varphi_k}_{g_0}^2 \psi^2\,dV_{g_0} \leq \frac{C}{k^2}\int_{\Sph^n} \big(\abs{\nabla \psi}_{g_0}^2+\psi^2\big)\,dV_{g_0}.
\end{equation*}
\end{prop}

This proposition will follow from the following two-dimensional Hardy estimate.

\begin{lem}[Critical Hardy inequality on a disk]\label{lemma.Hardy}
Let $R>0$ and denote by $B_R(0)\subset \R^2$ the two-dimensional disk of radius $R$. Then, for every $h\in W^{1,2}_0(B_R(0))$, we have the estimate
\begin{equation*}
	\int_{B_R(0)} \frac{h^2}{\abs{x}^2\log^2(R/\abs{x})}\,dx \leq 4\int_{B_R(0)}\abs{\partial_r h}^2\,dx.
\end{equation*}
\end{lem}

\begin{proof}[Proof of Lemma \ref{lemma.Hardy}]
The result is well-known, see for example \cite[Theorem 5]{MOW} for a slightly more general version. We nevertheless give a proof here to make the exposition more self-contained, and because it is instructive to see how the logarithmic weight arises naturally from transforming the two-dimensional radial Dirichlet energy into a one-dimensional energy.

By density, it suffices to consider $h\in C_c^\infty(B_R(0))$. Recall the classical one-dimensional Hardy inequality which states that for an absolutely continuous function $H(t)$ with $H(0)=0$, we have
\begin{equation*}
	\int_0^\infty \left| \frac{H(t)}{t}\right|^p dt 
	\leq \left(\frac{p}{p-1}\right)^p \int_0^\infty \abs{\partial_t H(t)}^p \, dt
\end{equation*}
for any $p>1$. Thus, writing $x=r\vartheta$ (where $\vartheta \in \Sph^1$) and substituting
\begin{equation*}
	t=\log(R/r),\qquad H_\vartheta(t)=h(Re^{-t}\vartheta),
\end{equation*}
we find (with $p=2$)
\begin{equation*}
	\int_{B_R(0)} \frac{h^2}{\abs{x}^2\log^2(R/\abs{x})}\,dx 
	=\int_{\mathbb S^1}\int_0^\infty \frac{H_\vartheta(t)^2}{t^2}\,dt\,d\vartheta
	\leq 4\int_{\mathbb S^1}\int_0^\infty \abs{\partial_t H_\vartheta(t)}^2\,dt\,d\vartheta
	=4\int_{B_R(0)} \abs{\partial_r h}^2\,dx.
\end{equation*}
This proves the lemma.
\end{proof}

\begin{proof}[Proof of Proposition \ref{prop.gradient}]
Pick a smooth cutoff function $\eta=\eta(r)$ such that $\eta=1$ on $\{r\leq 1/4\}$ and $\eta=0$ on $\{r\geq 3/8\}$. We first claim that, for every smooth function $\psi$ on $\Sph^n$,
\begin{equation}\label{eq.weighted-hardy}
	\int_{\{r<1/4\}} \frac{\psi^2}{r^2\log^2(1/r)}\,dV_{g_0} 
	\leq C\int_{\Sph^n} \left( \abs{\nabla \psi}_{g_0}^2+\psi^2\right)\,dV_{g_0}.
\end{equation}
Indeed, applying Lemma \ref{lemma.Hardy} with $R=1/2$ to $h = \eta \psi$ on each two-dimensional normal disk $B_{1/2}^\perp(y)$, and integrating over $y\in\Sigma$, we obtain
\begin{align*}
	\int_{\{r<1/4\}} \frac{\psi^2}{r^2\log^2(1/r)}\,dV_{g_0}
	&\leq C\int_\Sigma\int_{B_{1/2}^\perp(y)} \frac{(\eta \psi)^2}{r^2\log^2(1/(2r))} \,dx\,dV_\Sigma\\
	&\leq C\int_\Sigma\int_{B_{1/2}^\perp(y)} \abs{\partial_r(\eta \psi)}^2\,dx\,dV_\Sigma\\
	&\leq C\int_{\Sph^n} \bigl(\abs{\nabla \psi}_{g_0}^2+\psi^2\bigr)\,dV_{g_0}.
\end{align*}
Here we used the uniform comparability of $dV_{g_0}$ and $dx\,dV_\Sigma$ in tubular coordinates, together with $\log(1/r)\geq\log(1/(2r))$ for $r<1/2$. This establishes \eqref{eq.weighted-hardy}. On the other hand, by the definition of the conformal factors $\varphi_k$,
\begin{equation*}
	\abs{\nabla \varphi_k}_{g_0} \leq \frac{c_n\norm{\chi'}_\infty}{k\,r\log(1/r)}
\end{equation*}
on $\operatorname{supp}(\nabla \varphi_k) \subseteq \{r \leq r_k\}\subseteq \{r<1/4\}$, for all sufficiently large $k$. Estimate \eqref{eq.weighted-hardy} thus yields
\begin{equation*}
	\int_{\Sph^n}\abs{\nabla \varphi_k}_{g_0}^2 \psi^2\,dV_{g_0} 
	\leq \frac{c_n^2\norm{\chi'}_\infty^2}{k^2} \int_{\{r<1/4\}} \frac{\psi^2}{r^2\log^2(1/r)}\,dV_{g_0}
	\leq \frac{C}{k^2} \int_{\Sph^n} \left( \abs{\nabla \psi}_{g_0}^2+\psi^2\right)\,dV_{g_0},
\end{equation*}
as claimed.
\end{proof}

\begin{prop}[Convergence of the entropies]\label{prop.conv-entropies}
The entropies of $g_k=e^{2\varphi_k}g_0$ on $\Sph^n$ converge to the entropy of the round metric $g_0$, i.e.
\begin{equation}
\nu(g_k) \to \nu_n = \nu(g_0), \qquad (k\to \infty).
\end{equation}
\end{prop}

\begin{proof}
Note that 
\begin{equation*}
	\int_{\Sph^n} v^2 \,dV_{g_k}=\int_{\Sph^n} v^2 e^{n\varphi_k} dV_{g_0} = \int_{\Sph^n} (e^{n\varphi_k/2}v)^2 dV_{g_0}
\end{equation*}
and hence it is often useful to switch between $v$ and 
\begin{equation}\label{eq.defw}
	w := e^{n\varphi_k/2}v. 
\end{equation}
In particular, $\norm{v}_{L^2(g_k)} = 1$ if and only if $\norm{w}_{L^2(g_0)} = 1$. This lets us rewrite the two integral parts of the $\W$-entropy from \eqref{eq.W-decomp2} in terms of the fixed round metric $g_0$. Using the conformal transformation formula for the scalar curvature 
\begin{equation*}
	R_{g_k} = e^{-2\varphi_k} \big(R_{g_0} - 2(n-1)\Delta_{g_0}\varphi_k 
	- (n-1)(n-2)\abs{\nabla \varphi_k}^2_{g_0}\big)
\end{equation*}
and integrating by parts, we obtain
\begin{equation}\label{eq.conformal-Q}
\begin{aligned}
	Q_{g_k}(v) &=\int_{\Sph^n} (4\abs{\nabla v}_{g_k}^2+R_{g_k} v^2) dV_{g_k}\\
	&=\int_{\Sph^n} e^{-2\varphi_k} \big( 4\abs{\nabla w}_{g_0}^2+R_{g_0} w^2 
	-4w\langle\nabla w,\nabla\varphi_k\rangle_{g_0} 
	-(n-2)w^2 \abs{\nabla\varphi_k}_{g_0}^2 \big) \,dV_{g_0}\\
	&= Q_{g_0}(w) 
	+ \int_{\Sph^n} \big(e^{-2\varphi_k}-1\big) 
	\big( 4\abs{\nabla w}_{g_0}^2+R_{g_0} w^2\big) \, dV_{g_0}\\ 
	&\quad -  \int_{\Sph^n} e^{-2\varphi_k} \big(4w\langle\nabla w,\nabla\varphi_k\rangle_{g_0} 
	+(n-2)w^2 \abs{\nabla\varphi_k}_{g_0}^2 \big) \,dV_{g_0}.
\end{aligned}
\end{equation}
Moreover, the Nash entropy type term satisfies
\begin{equation}\label{eq.conformal-N}
\begin{aligned}
	\Nash_{g_k}(v^2) &=\int_{\Sph^n} v^2\log v^2 \, dV_{g_k}
	=\int_{\Sph^n} w^2\log(e^{-n\varphi_k}w^2) \, dV_{g_0}\\
	&= \int_{\Sph^n} w^2\log w^2 \,dV_{g_0} - n \int_{\Sph^n} \varphi_k w^2 \,dV_{g_0}\\
	&= \Nash_{g_0}(w^2) - n \int_{\Sph^n} \varphi_k w^2 \,dV_{g_0}\\
\end{aligned}
\end{equation}
\emph{Step 1: upper bound for $\nu(g_k)$.} 

In this step, we prove 
\begin{equation}\label{eq.limsup-nu}
	\limsup_{k\to\infty}\nu(g_k)\leq\nu_n.
\end{equation}
Of course, if one could prove that the round sphere is the \emph{global maximiser} of the entropy, then this step would not be necessary. Recall from the previous section that, denoting by $V_n=\Vol(\Sph^n,g_0)$ the volume of the round sphere, the function $\bar{w}$ and scale $\bar{\tau}$ realising the infimum $\nu_n = \nu(g_0) = \W(g_0,\bar{w},\bar{\tau})$ are 
\begin{equation*}
	\bar{w} \equiv V_n^{-1/2},\qquad\text{and}\quad \bar{\tau}= \frac{1}{2(n-1)}.
\end{equation*}
A simple test function for $\nu(g_k)$ is therefore the corresponding $\bar{v}_k=V_n^{-1/2}e^{-n\varphi_k/2}$. By \eqref{eq.conformal-Q} and \eqref{eq.conformal-N}, this test function satisfies
\begin{equation*}
	Q_{g_k}(\bar{v}_k) = Q_{g_0}(\bar{w}) 
	+ \frac1{V_n} \int_{\Sph^n} \big( e^{-2\varphi_k} -1\big) R_{g_0} \,dV_{g_0}
	- \frac1{V_n} \int_{\Sph^n} e^{-2\varphi_k}(n-2)\abs{\nabla\varphi_k}_{g_0}^2 \,dV_{g_0}
\end{equation*}
as well as
\begin{equation*}
	\Nash_{g_k}(\bar{v}_k^2) = \Nash_{g_0}(\bar{w}^2) -\frac{n}{V_n} \int_{\Sph^n}\varphi_k\,dV_{g_0}.
\end{equation*}
Using \eqref{eq.W-decomp1} and $\bar{\tau}R_{g_0}=n/2$, we obtain
\begin{equation*}
	\W(g_k,\bar{v}_k,\bar{\tau}) = \nu_n 
	+ \frac{n}{2V_n}\int_{\Sph^n} \big( e^{-2\varphi_k} -1+2\varphi_k \big)\,dV_{g_0}
	- \frac{n-2}{2(n-1)V_n} \int_{\Sph^n}e^{-2\varphi_k} \abs{\nabla\varphi_k}_{g_0}^2 \,dV_{g_0}.
\end{equation*}
Since $-c_n\leq\varphi_k\leq0$ and $\varphi_k$ vanishes outside the tube $\{r<r_k\}$, the first integral is bounded above by
\begin{equation*}
	C\Vol_{g_0}(\{r<r_k\})\leq Cr_k^2.
\end{equation*}
The second integral is nonnegative and appears with a negative sign, it can therefore simply be dropped. We conclude that
\begin{equation}\label{eq.entropy-upper}
	\nu(g_k)\leq \W(g_k,\bar{v}_k,\bar{\tau}) \leq\nu_n+Cr_k^2,
\end{equation}
and in particular \eqref{eq.limsup-nu} follows. We remark that by exploiting the negative term rather than dropping it, one can actually show that $\nu(g_k)<\nu_n$ for sufficiently large $k$. We leave the details of this fact to the reader.

\emph{Step 2: lower bound for $\W_{g_k}(v,\tau)$.} 

Here we prove a rough lower bound for the $\W$-entropy for arbitrary test functions $v$ and scales $\tau$. Since $e^{-2\varphi_k} - 1 \geq 0$
and $R_{g_0}=n(n-1)>0$, \eqref{eq.conformal-Q} yields
\begin{equation*}
	Q_{g_k}(v) \geq Q_{g_0}(w) 
	-  \int_{\Sph^n} e^{-2\varphi_k} \big(4w\langle\nabla w,\nabla\varphi_k\rangle_{g_0} 
	+(n-2)w^2 \abs{\nabla\varphi_k}_{g_0}^2 \big) \,dV_{g_0}.
\end{equation*}
Note that because $e^{-2\varphi_k} \leq e^{2c_n}$, the error term satisfies
\begin{align*}
	& \left| \int_{\Sph^n} e^{-2\varphi_k} \big( 4w\langle\nabla w,\nabla\varphi_k\rangle_{g_0} 
	+(n-2)w^2 \abs{\nabla\varphi_k}_{g_0}^2 \big) \,dV_{g_0} \right| \\
	&\qquad \leq C \left(\int_{\Sph^n} \abs{\nabla w}_{g_0}^2 \, dV_{g_0} \right)^{1/2}
	\left(\int_{\Sph^n} w^2\abs{\nabla \varphi_k}_{g_0}^2 \,dV_{g_0} \right)^{1/2} 
	+ C \int_{\Sph^n} w^2\abs{\nabla \varphi_k}_{g_0}^2 \,dV_{g_0}\\
	&\qquad \leq \frac{C}{k} Q_{g_0}(w) + \frac{C}{k^2} Q_{g_0}(w)
\end{align*}
where the last estimate follows by Proposition \ref{prop.gradient}, for sufficiently large $k$. In particular, we find
\begin{equation}\label{eq.Q-comparison}
	Q_{g_k}(v)\geq \Big(1-\frac{C}{k}\Big) Q_{g_0}(w)
\end{equation}
for sufficiently large $k$, where the constant $C$ is independent of $k$. Write $\alpha_k=1-C/k$ with $\alpha_k \to 1$ as $k \to \infty$ and, in particular, $\alpha_k>0$ for sufficiently large $k$. Combining \eqref{eq.W-decomp1} with \eqref{eq.conformal-N} and \eqref{eq.Q-comparison}, we find that for every normalised smooth test function $v$ and every scale $\tau>0$,
\begin{equation}\label{eq.W-comparison}
\begin{aligned}
	\W_{g_k}(v,\tau)
	&\geq \alpha_k\tau Q_{g_0}(w) -\Nash_{g_0}(w^2) -n-\frac n2\log(4\pi\tau) 
	+n\int_{\Sph^n} \varphi_kw^2\,dV_{g_0}\\
	&= \W_{g_0}(w,\alpha_k\tau) +\frac n2\log\alpha_k +n\int_{\Sph^n}\varphi_kw^2\,dV_{g_0}.
\end{aligned}
\end{equation}
Since $\varphi_k\geq -c_n$ and $\W_{g_0}(w,s)\geq\nu_n$ for all $s>0$, this implies
\begin{equation*}
	\W_{g_k}(v,\tau) \geq \nu_n+\frac{n}{2}\log\alpha_k - n c_n.
\end{equation*}
In particular, $\nu(g_k)>-\infty$ for all sufficiently large $k$. 

\emph{Step 3: lower bound for $\nu(g_k)$.} 

In this step, we prove 
\begin{equation}\label{eq.liminf-nu}
	\liminf_{k\to\infty}\nu(g_k)\geq\nu_n
\end{equation}
which together with \eqref{eq.limsup-nu} finishes the proof of the proposition. To obtain \eqref{eq.liminf-nu}, we need to control the (almost) minimising sequence of scales and test functions. To this end, for each sufficiently large $k$, choose $\tau_k>0$ and a smooth normalised test function $v_k$ such that
\begin{equation}
	\nu(g_k)\leq \W_{g_k}(v_k,\tau_k) \leq\nu(g_k)+\frac{1}{k}.
\label{eq.almost-minimisers}
\end{equation}
Consider the corresponding $w_k=e^{n\varphi_k/2}v_k$, satisfying $\norm{w_k}_{L^2(g_0)}=1$. By the lower bound \eqref{eq.W-comparison}, the almost-minimising assumption \eqref{eq.almost-minimisers}, and the upper bound \eqref{eq.entropy-upper},
\begin{equation}\label{eq.scale-bound}
	\mu(g_0,\alpha_k\tau_k) +\frac{n}{2} \log\alpha_k-nc_n \leq \W_{g_k}(v_k,\tau_k) 
	\leq \nu_n+\frac{1}{k}+Cr_k^2.
\end{equation}
Recall from Perelman \cite[Section 3.1]{Per} that
\begin{equation*}
	\lim_{s\searrow0}\mu(g_0,s)=0, \qquad \lim_{s\nearrow\infty}\mu(g_0,s)=+\infty.
\end{equation*}
By our choice $0<c_n<-\nu_n/n$, we have $-nc_n>\nu_n$. Consequently, \eqref{eq.scale-bound} excludes both $\tau_k\to0$ and $\tau_k\to\infty$ along any subsequence. Thus there are constants $0<a<b<\infty$, independent of $k$, such that
\begin{equation}\label{eq.scale-compactness}
	a\leq\tau_k\leq b
\end{equation}
for all sufficiently large $k$.

Next, we claim that for some uniform constant $C$, we have
\begin{equation}\label{eq.Q-bound}
	Q_{g_0}(w_k)\leq C.
\end{equation}
Indeed, the definition of $\nu_n$ gives $\W_{g_0}(w_k,a/4)\geq\nu_n$, which can be rearranged as
\begin{equation*}
	\Nash_{g_0}(w_k^2) \leq \frac{a}{4} Q_{g_0}(w_k) -\nu_n-n-\frac n2\log(\pi a).
\end{equation*}
Substituting this into the first line of \eqref{eq.W-comparison}, and using $\varphi_k\geq-c_n$, yields
\begin{equation*}
	\W_{g_k}(v_k,\tau_k) - \nu_n-\frac n2\log\frac{a}{4\tau_k} +nc_n
	\geq \left(\alpha_k\tau_k-\frac{a}{4} \right) Q_{g_0}(w_k).
\end{equation*}
For large enough $k$, we have $\alpha_k\geq 1/2$ and together with \eqref{eq.scale-compactness}, we conclude $\alpha_k\tau_k-\frac{a}{4} \geq \frac{a}{4}$. This gives the desired upper bound \eqref{eq.Q-bound} due to the uniform upper bound for the terms on the left hand side.

We now estimate
\begin{equation*}
	\left| \int_{\Sph^n}\varphi_kw_k^2\,dV_{g_0} \right| \leq c_n\int_{\{r<r_k\}}w_k^2\,dV_{g_0} 
	\leq c_n r_k^2\log^2(1/r_k) \int_{\{r<r_k\}} \frac{w_k^2}{r^2\log^2(1/r)}\,dV_{g_0}.
\end{equation*}
Here we used the monotonicity of $r^2\log^2(1/r)$ for sufficiently small $r$. By weighted Hardy estimate
\eqref{eq.weighted-hardy}, and using once again $R_{g_0}>0$, we therefore find
\begin{equation*}
	\left| \int_{\Sph^n}\varphi_kw_k^2\,dV_{g_0} \right| \leq Cr_k^2\log^2(1/r_k)Q_{g_0}(w_k).
\end{equation*}
By the uniform bound \eqref{eq.Q-bound} and $r_k^2\log^2(1/r_k)\to0$, this implies
\begin{equation*}
	\int_{\Sph^n}\varphi_k w_k^2\,dV_{g_0} \to 0, \qquad k\to \infty.
\end{equation*}
Applying \eqref{eq.W-comparison} once more gives
\begin{equation*}
	\nu(g_k)+\frac1k \geq\W_{g_k}(v_k,\tau_k) \geq\nu_n+\frac{n}{2} \log\alpha_k
	+n\int_{\Sph^n}\varphi_kw_k^2\,dV_{g_0} \to \nu_n,
\end{equation*}
which implies \eqref{eq.liminf-nu} and thus finishes the proof.
\end{proof}

We can now prove our instability theorem.

\begin{proof}[Proof of Theorem \ref{thm.main2}]
Let $d_k$ and $d_0$ denote the distance functions associated with $g_k=e^{2\varphi_k}g_0$ and $g_0$, respectively. Since $-c_n\leq\varphi_k\leq 0$, we have $e^{-2c_n}g_0\leq g_k\leq g_0$, and hence
\begin{equation}\label{eq.distance-comparison}
	e^{-c_n}d_0(x,y)\leq d_k(x,y)\leq d_0(x,y)
\end{equation}
for all $x,y\in\Sph^n$. Moreover, since $\varphi_k=-c_n$ on $\Sigma$ and $\Sigma$ is totally geodesic,
\begin{equation}\label{eq.distance-on-sigma}
	d_k(x,y)=e^{-c_n}d_0(x,y), \qquad x,y\in\Sigma.
\end{equation}

We claim that $(\Sph^n,g_k)$ remain a definite positive Gromov--Hausdorff distance away from $(\Sph^n,g_0)$. One possible way to establish this is to pass to a limit of the distance functions. Indeed, by \eqref{eq.distance-comparison} and the Arzel\`{a}--Ascoli theorem, a subsequence of $d_k$ converges uniformly to a metric $d_\infty$. Since the metrics $g_k$ agree with $g_0$ outside smaller and smaller tubular neighbourhoods of $\Sigma$, the limiting distance $d_\infty$ agrees \emph{locally} with $d_0$ on $\Sph^n\setminus\Sigma$. If $(\Sph^n,d_\infty)$ were isometric to $(\Sph^n,d_0)$, such an isometry would restrict to a local isometry of the round metric on $\Sph^n\setminus\Sigma$. By connectedness of $\Sph^n \setminus \Sigma$ and the rigidity of local isometries of the round sphere, this restriction would agree with a global round isometry. Continuity would then imply $d_\infty=d_0$ \emph{everywhere}, contradicting \eqref{eq.distance-on-sigma}. The disadvantage of this proof is that a contradiction argument only provides a qualitative statement.

Instead, we give a direct \emph{quantitative} argument here, using M\'{e}moli's lower bound for the Gromov--Hausdorff distance in terms of the radii of compact metric spaces \cite[Theorem 3.4]{Mem}. Recall that the radius of a compact metric space $(X,d)$ is defined by
\begin{equation*}
	\operatorname{rad}(X,d) =\inf_{x\in X}\sup_{y\in X}d(x,y).
\end{equation*}
For any two compact metric spaces $(X,d_X)$ and $(Y,d_Y)$, M\'{e}moli proves that
\begin{equation}\label{eq.GH-radius}
	\dGH(X,Y)\geq \frac12\abs{ \operatorname{rad}(X,d_X) -\operatorname{rad}(Y,d_Y)}.
\end{equation}

For the unit round sphere, every point has an antipodal point at distance $\pi$, so $\operatorname{rad}(\Sph^n,d_0)=\pi$. We now show that the radius of $(\Sph^n,d_k)$ is uniformly smaller. Fix a point $p\in\Sigma$. Since $\Sigma$ is a totally geodesic subsphere, its antipodal point $-p$ also belongs to $\Sigma$, and
\eqref{eq.distance-on-sigma} gives
\begin{equation*}
	d_k(p,-p)=e^{-c_n}\pi.
\end{equation*}
For any $q\in\Sph^n$, set $t=d_0(p,q)$. Since $d_0(-p,q)=\pi-t$, the triangle inequality and \eqref{eq.distance-comparison} imply
\begin{equation*}
	d_k(p,q) \leq \min \{ d_0(p,q), \, d_k(p,-p)+d_0(-p,q) \} =\min \{ t, \, e^{-c_n}\pi+\pi-t \}.
\end{equation*}
As the minimum of two numbers is always smaller or equal to their average, this yields
\begin{equation*}
	d_k(p,q) \leq \frac{\pi}{2}(1+e^{-c_n}).
\end{equation*}
Thus
\begin{equation*}
	\operatorname{rad}(\Sph^n,d_k) \leq \frac{\pi}{2}(1+e^{-c_n})
\end{equation*}
and \eqref{eq.GH-radius} yields the explicit lower bound
\begin{equation}\label{eq.GH-lower-bound}
    \dGH \big((\Sph^n,g_k),(\Sph^n,g_0)\big) \geq \frac{\pi}{4}(1-e^{-c_n})>0.
\end{equation}

We have already proved $\nu(g_k) \to \nu_n$ as $k \to \infty$ and $-c_n\leq\varphi_k\leq0$ gives uniform bi-Lipschitz equivalence. The metrics constructed have volumes approaching that of the round sphere, as
\begin{equation*}
	\abs{\Vol(\Sph^n,g_k)-V_n} \leq C\Vol_{g_0}(\{r<r_k\}) \to 0.
\end{equation*}
In order to have metrics that have \emph{precisely} the volume $V_n$, we can simply set
\begin{equation*}
	a_k=\left(\frac{V_n}{\Vol(\Sph^n,g_k)}\right)^{2/n} \to 1,\qquad \widehat g_k=a_kg_k.
\end{equation*}
Then $\Vol(\Sph^n,\widehat g_k)=V_n$, and $\nu(\widehat g_k)=\nu(g_k)$ by scaling-invariance of the entropy. Conformality and uniform bi-Lipschitz equivalence to the round metric persist. Finally,
\begin{equation*}
	\dGH \big((\Sph^n,\widehat g_k),(\Sph^n,g_k)\big) 
	\leq \frac{1}{2}\abs{ \sqrt{a_k}-1}\operatorname{diam}(\Sph^n,g_k) \to 0,
\end{equation*}
Consequently, the definite separation from the round sphere survives. More precisely, \eqref{eq.GH-lower-bound} implies
\begin{equation*}
	\liminf_{k\to\infty} \dGH \big((\Sph^n,\widehat g_k),(\Sph^n,g_0)\big) \geq \frac{\pi}{4}(1-e^{-c_n})>0,
\end{equation*}
which completes the proof.
\end{proof}

\begin{rem}
The preceding argument also shows that the metrics $g_k$ remain uniformly separated from round spheres of \emph{any} radius. Indeed, let $p\in\Sph^n$ be a point with $d_0(p,\Sigma)=\pi/2$. Since $g_k=g_0$ outside $\{r<r_k\}$, any curve joining $p$ to its antipode $-p$ has $g_k$-length at least $\pi-2r_k$. Consequently,
\begin{equation*}
	\pi-2r_k \leq \operatorname{diam}(\Sph^n,g_k) \leq \pi,
\end{equation*}
so $\operatorname{diam}(\Sph^n,g_k)\to\pi$. On the other hand, we have shown that
\begin{equation*}
	\operatorname{rad}(\Sph^n,g_k) \leq \frac{\pi}{2}(1+e^{-c_n}).
\end{equation*}
Since every round sphere has equal radius and diameter, M\'{e}moli's estimate \cite[Theorem 3.4]{Mem} gives, for every radius $\varrho>0$,
\begin{equation*}
	\dGH\big((\Sph^n,g_k),(\Sph^n,\varrho^2g_0)\big) 
	\geq \frac{1}{4} \big( \operatorname{diam}(\Sph^n,g_k)  -\operatorname{rad}(\Sph^n,g_k)\big) 
	\geq \frac{\pi}{8}(1-e^{-c_n})-\frac{r_k}{2}.
\end{equation*}
The same conclusion holds for the volume-normalised metrics $\widehat g_k$.
\end{rem}

\section{Gromov--Hausdorff stability of $\Sph^2$}\label{sec.stability}

In this section, we prove the Gromov--Hausdorff stability theorem in dimension two. The topological part will follow from the following more general result.

\begin{prop}[Simple entropy estimate based on topology]\label{prop.entropy-topology}
Let $(M^2,g)$ be a smooth closed Riemannian surface. If $\chi(M)\leq 0$, then $\nu(g)=-\infty$. If $\chi(M)>0$, then $\nu(g)\leq\log\chi(M)-1$. In particular,
\begin{align*}
	M\cong\mathbb S^2 &\implies \nu(g)\leq\log2-1,\\
	M\cong\mathbb{RP}^2 &\implies \nu(g)\leq-1.
\end{align*}
Both upper bounds are sharp, since equality holds for the corresponding round metrics.
\end{prop}

\begin{proof}
Let $A=\Area(M,g)$ and use the constant test function $v=A^{-1/2}$ in Perelman's $\W$-functional. By Gauss–Bonnet,
\begin{equation*}
	\int_M R_g\,dV_g=4\pi\chi(M).
\end{equation*}
Consequently,
\begin{equation*}
	\W_g(v,\tau) =\frac{4\pi\chi(M)}{A}\tau +\log A-2-\log(4\pi\tau) 
	=\frac{4\pi\chi(M)}{A}\tau -\log\frac{4\pi\tau}{A}-2.
\end{equation*}
If $\chi(M)\leq0$, letting $\tau\to\infty$ shows that $\nu(g)=-\infty$. If instead $\chi(M)>0$, then the right-hand side is minimised at $\tau_{min}=\frac{A}{4\pi\chi(M)}$ with value $\log\chi(M)-1$ and hence
\begin{equation*}
	\nu(g) \leq \log\chi(M)-1.
\end{equation*}
The sharpness of the bounds follows from the calculations for round spheres and their quotients in Section \ref{sec.preliminaries}.
\end{proof}

In addition, we will also use the following spectral rigidity result of Xu for Riemannian metrics $g_k$ on the two-sphere $\Sph^2$, which we can reformulate as follows.

\begin{prop}[Theorem 1.8 of \cite{Xu}]\label{prop.Xu}
Let $\beta>\frac{1}{4}$. For any sequence of metrics $h_k$ on $\Sph^2$ with $\lambda_1(-\Delta_{h_k}+\beta K_{h_k}) \geq \beta$ and $\Area(\Sph^2,h_k)\to 4\pi$, there exist diffeomorphisms $\Psi_k: \Sph^2 \to \Sph^2$ with
\begin{equation*}
	\max_{x,y\in \Sph^2} \abs{ d_{g_0}(\Psi_k(x),\Psi_k(y)) - d_{h_k}(x,y) } =\eps_k \to 0.
\end{equation*}
\end{prop} 

In this proposition $K_{h_k} = \frac{1}{2}R_{h_k}$ denotes the Gauss curvature. We will apply this result with $\beta=\frac{1}{2}>\frac{1}{4}$. As we will see, the first eigenvalue of 
\begin{equation*}
	-\Delta_{h_k}+\beta K_{h_k} = -\Delta_{h_k} + \frac{1}{4} R_{h_k}
\end{equation*}
can be controlled by the $\nu$-entropy.

We are now ready to prove the Gromov--Hausdorff stability theorem.

\begin{proof}[Proof of Theorem \ref{thm.main1}]
Let $(M_k,g_k)_{k \in \N}$ be a sequence of smooth, closed Riemannian surfaces satisfying $\nu(g_k)\to \nu_2 = \log 2-1 \approx -0.30685$ as well as the normalisation
\begin{equation*}
	\Area(M_k,g_k)=4\pi.
\end{equation*}
For sufficiently large $k$, $\nu(g_k) > -1$ and therefore by Proposition \ref{prop.entropy-topology}, $M_k$ is diffeomorphic to $\Sph^2$.

Next, set $\lambda_k=\lambda_1(-4\Delta_{g_k}+R_{g_k})$, and choose a positive first eigenfunction $v_k$ normalised by $\int v_k^2\,dV_{g_k}=1$. We claim that $\Nash_{g_k}(v_k^2)\geq-\log(4\pi)$.
Indeed, recall that Jensen's inequality
\begin{equation*}
	\Phi\left(\int_{\Sph^2}  h\,d\mu\right)\leq\int_{\Sph^2} \Phi(h)\,d\mu,
\end{equation*}
holds for any probability measure $\mu$ and any convex function $\Phi$. Applying it to $\Phi(t)=t\log t$, $h=4\pi v_k^2$ and $d\mu=dV_{g_k}/(4\pi)$ and using $\int_{\Sph^2} h\,d\mu=1$, we get 
\begin{equation*}
	0=1\log 1=\Phi(1) \leq\int_{\Sph^2}(4\pi v_k^2)\log(4\pi v_k^2) \frac{dV_{g_k}}{4\pi} =\log(4\pi)+\Nash_{g_k}(v_k^2).
\end{equation*}

Since $Q_{g_k}(v_k)=\lambda_k$, testing the entropy with $v_k$ yields, for every $\tau>0$, 
\begin{equation*}
	\nu(g_k)\leq\W_{g_k}(v_k,\tau) \leq\tau\lambda_k-\log\tau-2.
\end{equation*}
If $\lambda_k\leq0$, letting $\tau\to\infty$ contradicts the finite lower bound on $\nu(g_k)$, thus $\lambda_k>0$. In fact, choosing $\tau=1/\lambda_k$ and $\delta_k=\nu_2-\nu(g_k)$ gives
\begin{equation*}
	\log 2-1-\delta_k= \nu(g_k)\leq\log\lambda_k-1, 
\end{equation*}
and hence
\begin{equation*}
	\lambda_k\geq 2e^{-\delta_k}.
\end{equation*}
Moreover, by Gauss--Bonnet
\begin{equation*}
	\lambda_k \leq\frac{1}{4\pi}\int_{\Sph^2} R_{g_k}\,dV_{g_k} =\chi(M_k)=2.
\end{equation*}
Consequently,
\begin{equation}
	2e^{-\delta_k} \leq \lambda_k \leq 2.
\end{equation}
Set $a_k=\lambda_k/2$ and define the metrics $h_k=a_k g_k$, so $a_k \to 1$. By scaling, we have
\begin{equation*}
	\lambda_1\Big( -\Delta_{h_k}+\frac{1}{4}R_{h_k} \Big)=\frac{\lambda_k}{4a_k}
	=\frac{1}{2},\qquad\Vol(\Sph^2,h_k)=4\pi a_k\to 4\pi.
\end{equation*}
Applying Proposition \ref{prop.Xu} to these metrics gives diffeomorphisms
$\Psi_k:\Sph^2 \to\Sph^2$ such that 	
\begin{equation*}
	\max_{x,y\in \Sph^2} \abs{ d_{g_0}(\Psi_k(x),\Psi_k(y)) - d_{h_k}(x,y) } =\eps_k \to 0.
\end{equation*}
Finally, $d_{g_k}=a_k^{-1/2}d_{h_k}$ and $\diam(\Sph^2,g_0)=\pi$, hence the triangle inequality gives 
\begin{equation*}
	\max_{x,y\in \Sph^2} \abs{ d_{g_0}(\Psi_k(x),\Psi_k(y)) - d_{g_k}(x,y) } \leq a_k^{-1/2}\eps_k 
	+\pi \abs{a_k^{-1/2}-1} \to 0.
\end{equation*}
This proves the Gromov--Hausdorff convergence.
\end{proof}

\section{Ricci flow stability of $\Sph^3$}\label{sec.dynamical}

In this section, we prove Theorem~\ref{thm.main3}. We first collect two useful lemmas. The first one states that the $\nu_c$-entropy is upper semicontinuous under pointed smooth convergence.

\begin{lem}[Entropy under pointed convergence]\label{lemma.entropy-pointed}
Let $(M_k,g_k,p_k)$ be a sequence of complete pointed Riemannian $n$-manifolds converging in the pointed smooth Cheeger--Gromov sense to a complete pointed manifold $(N,h,p_\infty)$. Then
\begin{equation*}
	\limsup_{k\to\infty}\nu_c(M_k,g_k)\leq\nu_c(N,h).
\end{equation*}
\end{lem}

\begin{proof}
This follows by a standard test-function argument, cf.~\cite[Lemma~6.28]{CCG} for the corresponding statement for $\mu$. We briefly sketch the argument here. Fix $\tau>0$ and a nonnegative function $v\in C_c^\infty(N)$ with $\int_Nv^2\,dV_h=1$. By the smooth Cheeger--Gromov convergence, there exists an exhaustion of $N$ by open sets $U_k$ and embeddings $\Phi_k: U_k \to M_k$ such that $\Phi_k^* g_k$ smoothly converges to $h$. In particular, for sufficiently large $k$, the support of $v$ lies in $U_k$, and hence we can push-forward $v$ to $M_k$ with $\Phi_k$, extend it by zero and normalise it so that we obtain a smooth test function $v_k$ with $\norm{v_k}_{L^2(g_k)}=1$. By smooth convergence, one obtains
\begin{equation*}
	\W_{g_k}(v_k,\tau)\to\W_h(v,\tau).
\end{equation*}
On the other hand $\nu_c(M_k,g_k)\leq\W_{g_k}(v_k,\tau)$. The claim then follows by taking the $\limsup$ in this inequality, followed by taking the infimum over all test functions $v$ and scales $\tau$.
\end{proof}

The second lemma is a consequence of Perelman's neck theorem \cite[Section 11.8]{Per} for noncompact three-dimensional $\kappa$-solutions. Recall that given $\kappa>0$, a three-dimensional $\kappa$-solution on $(-\infty,0]$ is a three-dimensional Ricci flow $(N^3,h(t))_{t\leq0}$ that is $\kappa$-noncollapsed at all scales, for which $\abs{\Rm}$ is uniformly bounded on each compact time interval in $(-\infty,0]$, and each time-slice is complete and has nonnegative curvature operator and positive scalar curvature.

\begin{lem}[Entropy bound for noncompact $\kappa$-solutions]
\label{lem.entropy-kappa}
Let $(N^3,h(t))_{t\leq0}$ be an oriented, nonflat, noncompact $\kappa$-solution. Then
\begin{equation}\label{eq.kappa-cylinder-bound}
\nu_c(N,h(0))\leq\nu_2=\log2-1.
\end{equation}
\end{lem}

\begin{proof}
By Perelman's neck theorem (see Kleiner-Lott \cite[Corollary~47.2]{KL} for a more detailed version), there are points $q_k\in N$ and constants $\lambda_k>0$ such that
\begin{equation*}
	(N,\lambda_k h(0),q_k) \to (\Sph^2\times\R,h_{\mathrm{cyl}},q_\infty)
\end{equation*}
smoothly in the pointed Cheeger--Gromov sense, where $h_{\mathrm{cyl}}=2g_0+dz^2$ and $g_0$ is the unit round metric on $\Sph^2$. Recall that in Section \ref{sec.preliminaries} we considered the metric $g_{\mathrm{cyl}}=g_0+dz^2$. Scale invariance, Lemma~\ref{lemma.entropy-pointed}, and \eqref{eq.cylinder-entropy} give
\begin{equation*}
	\nu_c(N,h(0)) \leq\nu_c(\Sph^2\times\R,h_{\mathrm{cyl}}) = \nu_c(\Sph^2\times\R,g_{\mathrm{cyl}})
	\leq \nu_2 = \log 2 -1.
\end{equation*}
This finishes the proof.
\end{proof}

We are now ready to prove our last main result.

\begin{proof}[Proof of Theorem~\ref{thm.main3}]
Let $(M,g)$ be a smooth, closed, connected Riemannian three-manifold satisfying the assumption
\begin{equation}\label{eq.high-entropy}
\nu(M,g)>\nu_2.
\end{equation}
We will treat the orientable and non-orientable cases separately.

\emph{Step 1: the orientable case.} 

Assume that $M$ is orientable. Let $\lambda_1(g)$ be the lowest eigenvalue of $-4\Delta_g+R_g$, and let $v$ be a corresponding positive $L^2(g)$-normalised eigenfunction. If $\lambda_1(g)\leq 0$, then, as $\tau \to \infty$, 
\begin{align*}
	\W_g(v,\tau) =\tau\lambda_1(g)-\Nash_g(v^2)-3-\frac{3}{2} \log (4\pi\tau) \to -\infty.
\end{align*}
contradicting \eqref{eq.high-entropy}. Thus $\lambda_1(g)>0$. For the maximal unnormalised Ricci flow $g(t)$, $t\in[0,T)$, Perelman's evolution inequality \cite{Per} gives
\begin{equation*}
	\frac{d}{dt}\lambda_1(g(t))\geq\frac{2}{3}\lambda_1(g(t))^2,
\end{equation*}
and we can conclude that $T\leq\frac{3}{2\lambda_1(g)}<\infty$. In other words, the entropy assumption \eqref{eq.high-entropy} implies that the Ricci flow starting at $(M,g)$ must develop a finite time singularity.

Now pick times $t_k\nearrow T$ and points $p_k\in M$ with the property that
\begin{equation*}
	Q_k:=|\Rm|_{g(t_k)}(p_k) =\max_{M\times[0,t_k]}\abs{\Rm}\to \infty,
\end{equation*}
and consider the rescaled flows
\begin{equation*}
	g_k(s)=Q_k g\big( t_k+\tfrac{s}{Q_k}\big), \qquad s\in [- Q_kt_k,0].
\end{equation*}
The rescaled curvatures satisfy $\abs{\Rm}_{g_k(s)}\leq1$ and $\abs{\Rm}_{g_k(0)}(p_k)=1$. Perelman's three-dimensional singularity analysis (using his no-local-collapsing theorem, Hamilton--Ivey pinching, and Hamilton's compactness theorem) therefore yields a subsequence converging smoothly in the pointed Cheeger--Gromov sense to a complete, nonflat, oriented ancient three-dimensional $\kappa$-solution $(N,h(s),p_\infty)$, $s\leq0$; see \cite{Per,KL}. Monotonicity and scale invariance of the entropy, as well as Lemma~\ref{lemma.entropy-pointed} then imply
\begin{equation}\label{eq.entropy-blowup}
	\nu_2 < \nu(g) \leq\limsup_{k\to\infty}\nu(g_k(0)) \leq\nu_c(N,h(0)).
\end{equation}
Lemma~\ref{lem.entropy-kappa} then forces $N$ to be compact. In particular, the convergence is in the global smooth Cheeger--Gromov sense, that is, the embeddings $\Phi_k: N \to (M,g_k)$ in the definition of pointed Cheeger--Gromov convergence eventually become global diffeomorphisms onto the approximating manifolds.

A compact three-dimensional $\kappa$-solution has strictly positive sectional curvature; see for example Kleiner-Lott \cite[Section~59]{KL} and hence strictly positive Ricci curvature. Consequently, global smooth convergence implies $\Ric_{g(t_k)}>0$ for sufficiently large $k$. Hamilton's theorem \cite{Ham82} then gives long-time existence and exponential convergence in $C^\infty$ of the volume-normalised Ricci flow starting from $(M,g(t_k))$ to a metric $g_\infty$ of constant positive sectional curvature. The rescaling relation between normalised and unnormalised Ricci flow then gives the same conclusion for the volume-normalised flow starting from $(M,g)$.

We conclude that $M$ is a spherical space form $\Sph^3/\Gamma$. Applying entropy monotonicity and Lemma~\ref{lemma.entropy-pointed} once again, together with \eqref{eq.quotient-entropy}, we obtain
\begin{equation}\label{eq.entropy-round-limit}
	\nu(g)\leq\nu(g_\infty)=\nu_3-\log|\Gamma|.
\end{equation}
If $\abs{\Gamma}\geq2$, this contradicts \eqref{eq.high-entropy}. Therefore, $\Gamma$ is trivial. This finishes the proof in the orientable case.

\emph{Step 2: the non-orientable case.} 

First note that for any finite Riemannian covering of closed manifolds $\pi:(\widetilde M,\widetilde g)\to(M,g)$,
\begin{equation}\label{eq.entropy-cover}
	\nu(M,g)\leq\nu(\widetilde M,\widetilde g).
\end{equation}
Indeed, for a smooth positive $L^2(\widetilde g)$-normalised test function $\widetilde v$ on $\widetilde M$, define $v$ by
\begin{equation*}
	v(y)^2=\sum_{x\in\pi^{-1}(y)}\widetilde v(x)^2.
\end{equation*}
Then $v$ is a smooth, positive, and $L^2(g)$-normalised function on $M$. By Cauchy--Schwarz and the elementary inequality
\begin{equation*}
	\Big(\sum_j a_j \Big)\log \Big(\sum_j a_j\Big) \geq\sum_j a_j\log a_j, \qquad a_j\geq 0,
\end{equation*}
we therefore have
\begin{equation*}
	\W_g(v,\tau)\leq\W_{\widetilde g}(\widetilde v,\tau).
\end{equation*}
Taking infima over $\widetilde v$ and $\tau>0$ proves \eqref{eq.entropy-cover}.

Now, if $M$ were non-orientable and satisfied \eqref{eq.high-entropy}, its connected orientable double cover would also have entropy above $\nu_2$ by \eqref{eq.entropy-cover}, and hence would be diffeomorphic to $\Sph^3$ by the orientable case from Step 1. But every manifold covered by $\Sph^3$ is orientable. In fact, a nontrivial deck transformation would be a fixed-point-free orientation-reversing diffeomorphism of $\Sph^3$. This is impossible by the Lefschetz fixed-point theorem, since its Lefschetz number is $1-(-1)=2$. This completes the proof of dynamical stability.

\emph{Step 3: maximality and rigidity.}

Maximality follows directly from the above steps. Indeed, for any smooth, closed Riemannian three-manifold, either $\nu(g)\leq\nu_2<\nu_3$, or the preceding arguments give $\nu(g)\leq\nu_3$.

To prove rigidity, suppose that $\nu(g)=\nu_3$. By monotonicity and maximality, the entropy remains constant along the Ricci flow, and Perelman's equality case \cite{Per} implies that $g$ is a compact gradient shrinking Ricci soliton. Consequently, its volume-normalised Ricci flow evolves only by diffeomorphisms. Since we have already shown that this flow converges smoothly to a round metric, $g$ itself must be round, up to scaling and isometry. This finishes the proof of the theorem.
\end{proof}

\makeatletter
\def\@listi{%
  \itemsep=0pt
  \parsep=1pt
  \topsep=1pt}
\makeatother
{\fontsize{10}{11}\selectfont

}
\printaddress

\end{document}